\documentclass[conference, onecolumn]{IEEEtran}

\usepackage{cite}
\usepackage{amsmath,amssymb,amsfonts}
\usepackage{algorithmic}
\usepackage{graphicx}
\usepackage{textcomp}
\usepackage{xcolor}
\usepackage[hidelinks]{hyperref}

\usepackage[numbers,sort&compress]{natbib}

\title{Vectorization of Persistence Diagrams for Topological Data Analysis in R and Python Using TDAvec Package}

\author{
\IEEEauthorblockN{Aleksei Luchinskii}
\IEEEauthorblockA{\textit{Bowling Green State University} \\
Bowling Green, OH, USA \\
aluchi@bgsu.edu}
\and 
\IEEEauthorblockN{Umar Islambekov}
\IEEEauthorblockA{\textit{Bowling Green State University} \\
Bowling Green, OH, USA \\
iumar@bgsu.edu}
}

\usepackage{graphicx}

\begin{document}

\maketitle{}
\begin{abstract}
Persistent homology is a widely-used tool in topological data analysis (TDA) for understanding the underlying shape of complex data. By constructing a filtration of simplicial complexes from data points, it captures topological features such as connected components, loops, and voids across multiple scales. These features are encoded in persistence diagrams (PDs), which provide a concise summary of the data's topological structure. However, the non-Hilbert nature of the space of PDs poses challenges for their direct use in machine learning applications. To address this, kernel methods and vectorization techniques have been developed to transform PDs into machine-learning-compatible formats. In this paper, we introduce a new software package designed to streamline the vectorization of PDs, offering an intuitive workflow and advanced functionalities. We demonstrate the necessity of the package through practical examples and provide a detailed discussion on its contributions to applied TDA. Definitions of all vectorization summaries used in the package are included in the appendix.
\end{abstract}

\section{Introduction}

Persistent homology has emerged as one of the most powerful tools in topological data analysis (TDA) for uncovering the geometric and topological structure of complex data \cite{Carlsson:2009, edelsbrunner2010computational, chazal2021introduction}. Within this framework, a data observation can take many forms—such as a point cloud in a metric space, an image, a graph, or a time series. The core idea of persistent homology is to construct a nested sequence (or filtration) of simplicial complexes indexed by a scale parameter, and to track how topological features appear and disappear across scales. These features correspond to “holes” of different dimensions—connected components, loops, voids, and their higher-dimensional analogues—captured through the algebraic machinery of homology. From a geometric perspective, simplicial complexes, composed of vertices, edges, triangles, and higher-dimensional simplices, provide a discrete representation that partially recovers the shape of the underlying space from sampled data \cite{nanda2013simplicial}.

The principal output of persistent homology is the persistence diagram (PD), a concise topological summary that encodes the lifetime of features across scales. A $k$-dimensional PD is defined as a multiset of points $D={(b_i,d_i)}_{i=1}^N$, where each point $(b_i,d_i)$ represents a feature of homological dimension $k$—with $k=0$ for connected components, $k=1$ for loops, and $k=2$ for voids. The birth coordinate $b_i$ indicates the scale at which a feature appears, and the death coordinate $d_i$ the scale at which it disappears.

In many practical applications, persistence diagrams serve as inputs to machine learning and statistical models. However, since PDs do not form a Hilbert space, they cannot be directly used with most standard algorithms that rely on vector or inner-product representations. To overcome this limitation, researchers have developed two major strategies: kernel methods and vectorization techniques \cite{chung2022persistence}. Kernel methods define similarity measures between pairs of PDs, while vectorization approaches transform diagrams into finite-dimensional feature vectors that can be readily processed by machine learning models. In the latter case, the procedure typically involves constructing a summary function from a PD—such as a persistence image, landscape, or entropy profile—and then discretizing it over a one- or two-dimensional grid of scale values. Both paradigms have been widely adopted and proven effective in applied TDA and data-driven research \cite{hensel2021survey}.

Despite the growing interest in topological data analysis, software support for vectorization of persistence diagrams remains fragmented across languages and packages, often differing in implementation details, computational efficiency, and interface design. To address this gap, we introduce TDAvec, a unified and efficient library for computing vector summaries of persistence diagrams in both R and Python. The package provides a consistent interface for several widely used vectorization techniques, allowing researchers to seamlessly integrate topological descriptors into machine learning workflows

The rest of the paper is organized as follows. In the next section we describe in more details our package and prove its necessity. In the next section a short overview of existing R and python packages is presented. In section \ref{sec:package-description} a short description of our package is given, followed by some usage examples. Last section is reserved for discussion. In the appendix you will find definitions of all used in our work vectorization summaries.

\section{Existing Models}

\subsection{R tools for TDA}

The computational tools for TDA in the \texttt{R} environment are provided through various packages\footnote{In this section, we only focus on \texttt{R} packages for TDA that are available on the CRAN repository.} such as \texttt{TDA} \cite{TDA}, \texttt{TDAstats} \cite{wadhwa2018tdastats}, \texttt{kernelTDA} \cite{kernelTDA}, \texttt{TDAmapper} \cite{TDAmapper}, \texttt{TDAkit} \cite{TDAkit}, \texttt{tdaunif} \cite{tdaunif}, \texttt{TDApplied} \cite{TDApplied} and \texttt{ripserr} \cite{ripserr} (see Table \ref{tab:TDA_packages} for an overview of these packages in terms of their scope and areas of focus). 

\newcommand{\cm}{$V$}

\newcommand{\rota}[1]{
  \rotatebox{90}{\texttt{#1}}
}

\begin{table*}
  \centering
  \begin{tabular}{|l|c|c|c|c|c|c|c|c||c|c|c|c|}
    \hline
    & \multicolumn{8}{c||}{R libraries} & \multicolumn{4}{c|}{Python} \\
    \hline
    & \rota{TDA} & \rota{TDApplied} & \rota{TDAstats} & \rota{TDAkit} & \rota{kernelTDA} &  \rota{tdaunif} &\rota{ripserr} & \rota{TDAmapper} & 
        \rota{Giotto-tda} &  \rota{Gudhi} &  \rota{Dionysus 2} & \rota{scikit-tda}\\
    \hline
    Sampling methods & \cm & & &\cm & & \cm& & &
      & & & \\
    \hline
    Density estimation & \cm & & & & & & & &
       & \cm & & \\
    \hline
    Alpha filtration & \cm & & & & & & & &
       & \cm & & \cm \\
    \hline
    Alpha shape filtration &\cm & & & & & & & &
       & & & \\
    \hline
    Vietoris-Rips filtration &\cm & &\cm &\cm & & & \cm& &
       & \cm & \cm & \cm \\
    \hline
    User-defined filtration &\cm & & & & & & & &
       & & & \\
    \hline
    Cubical complex & & & & & & & \cm & &
      \cm  & \cm & & \\
    \hline
    Wasserstein distance &\cm &\cm &\cm & & \cm & & & &
      \cm  & \cm & \cm & \cm \\
    \hline
    Plotting persistence diagrams &\cm &\cm &\cm &\cm & & & & &
      \cm  & \cm & \cm & \cm \\
    \hline
    Statistical methods &\cm &\cm &\cm &\cm & & & & &
       & & & \\
    \hline
    Vectorization methods &\cm & & &\cm & \cm &  & & &
       \cm & \cm & & \cm \\
    \hline
    Kernel methods & &\cm & &\cm & \cm & & & &
       & & & \\
    \hline
    Supervised learning methods & &\cm & & & \cm & & & &
       & & & \\
    \hline
    Clustering methods &\cm &\cm & & \cm& & & & &
      \cm  & \cm & & \\
    \hline
    HD data plots & & & & & & & & \cm &
     \cm   & \cm & \cm & \\
    \hline
    Dimension reduction & &\cm & & & & & & &
       & & & \cm \\
    \hline
  \end{tabular}
  \caption{\texttt{R} and \texttt{Python} packages for TDA}
  \label{tab:TDA_packages}
\end{table*}

The \texttt{TDA} package is the largest \texttt{R} package for TDA. The \texttt{TDA} package offers tools to compute PDs for commonly used types of filtrations such as \emph{Vietoris-Rips}, \emph{Alpha} and \emph{Alpha shape}. It also allows to construct more general sublevel set filtrations and compute the corresponding PDs. Moreover, the \texttt{TDA} package provides implementations to plot PDs and compute \emph{bottleneck and Wasserstein} distances between them. \texttt{TDAstats} offers a variety of tools for conducting statistical inference (such as hypothesis testing) on PDs. Compared to the \texttt{TDA} package, it computes PDs much faster for Vietoris-Rips filtrations based on the Ripser C++ library \cite{Bauer2021Ripser} and offers more aesthetic visualization of the diagrams using the \texttt{ggplot2} package \cite{Wickham2016ggplot2}. The \texttt{kernelTDA} package contains implementations of popular kernel-based methods for TDA such as \emph{geodesic Gaussian kernel}, \emph{geodesic Laplacian kernel}, \emph{persistence Fisher kernel} and \emph{persistence sliced Wasserstein kernel}. For computing the Wasserstein distance between a pair of PDs, unlike the \texttt{TDA} package, it uses an iterative procedure to reasonably approximate the exact distance which that leads to a considerable reduction in run-time cost. The \texttt{ripserr} package allows a fast computation of PDs for filtrations on Vietoris-Rips and cubical complexes using the Ripser C++ library. The \texttt{TDApplied} and \texttt{TDAkit} packages provides various tools to integrate topological features (PDs or their vector summaries) into machine and statistical learning settings. The \texttt{tdaunif} is a useful package if one needs to sample points from various manifolds such as a klein-bottle, an ellipse or a torus. \texttt{TDAmapper} offers tools to visualize high-dimensional data by constructing the so-called Mapper graphs that preserve its topological structure.

\subsection{Python tools for TDA}

Several \texttt{Python} libraries are available for TDA, such as \texttt{Giotto-tda} \cite{tauzin2021giotto}, \texttt{Ripser} \cite{christopher2018lean}, \texttt{Gudhi} \cite{rouvreau2020gudhi}, \texttt{Dionysus 2} \cite{Dionysus2}, \texttt{Scikit-tda} \cite{scikittda2019} which is a container library for such projects as  \texttt{Ripser}, \texttt{Persim} \cite{Persim},  \texttt{KeplerMapper} \cite{van2019kepler}, etc. \texttt{Giotto-tda} is a powerful library that integrates with the popular machine learning library \texttt{scikit-learn}, offering tools for persistent homology and visualizations of persistence diagrams. \texttt{Ripser} focuses on fast computation of Vietoris-Rips complexes, especially for large datasets. \texttt{Gudhi} provides a wide range of topological tools for simplicial complexes, persistent homology, and topological signatures. \texttt{Dionysus 2} offers fast computation of persistent homology and cohomology, with an emphasis on flexibility and efficiency. \texttt{Persim} focuses on tools for working with PDs. It contains implementations of commonly used vectorization and kernel methods for PDs. \texttt{KeplerMapper} implements the TDA Mapper algorithm to visualize high-dimensional data. \texttt{Scikit-tda} is another package that integrates with \texttt{scikit-learn}, simplifying the application of TDA to typical machine learning tasks. 

For a more comprehensive list of \texttt{Python} libraries for TDA and their functionality, we refer the readers to \cite{awesome-tda:2024}.

\section{Package Description}
\label{sec:package-description}

\subsection{Statement of need and Short Package Description}

As it was mentioned in the previous section, vectorization step is extremely important to include TDA in the ML pipeline. Up to now lots of different vectorization methods were proposed (see the Appendix below for short description of some of the existing ones). It turns out that performance accuracy of the ML algorithms depends strongly on the choice of the used vectorization, so it could be very interesting to be able to compare different approaches.

In order to do such a comparison it would be useful to have all considered vectorization methods implemented in one package in the uniform manner. Unfortunately, up to now such a package does not exist. Some mentioned in the Introduction libraries do have implementation of some of the vectorization (see table \ref{tab:vects} for quick comparison), but there is no package combining all of them together. More over, in R all the code behind the existing vector implementation is written using standard R factions, which may prove slow and inefficient in large-scale computations. In addition, the interfaces of various factions in different packages are not compatible with each other, which makes the comparison even more challenging.

\begin{table}[htbp]
  \centering
  \begin{tabular}{|l|c|c|c|c|c||c|c|c|}
    \hline
    & \multicolumn{5}{c||}{R libraries} & \multicolumn{3}{c|}{Python} \\
    \hline
    & \rota{TDA} & \rota{TDApplied} & \rota{TDAstats} & \rota{TDAkit} & \rota{kernelTDA} &   
        \rota{Giotto-tda} &  \rota{Gudhi} &  \rota{scikit-tda}\\
    \hline
    Bar Codes                   & \cm & \cm & \cm & \cm &     &     &     &      \\
    Persistence Landscape       & \cm &     &     & \cm &     & \cm & \cm &  \cm \\
    Persistence Silhouette      & \cm &     &     & \cm &     & \cm & \cm &      \\
    Persistent Entropy Summary  &     &     &     &     &     &     & \cm &      \\
    Betti Curve                 &     &     &     &     &     & \cm & \cm &      \\
    Euler characteristic curve  &     &     &     &     &     &     &     &      \\
    The normalized life curve   &     &     &     &     &     &     &     &      \\
    Persistence Surface         &     &     &     &     & \cm & \cm & \cm &  \cm \\
    Persistence Block           &     &     &     &     &     &     &     &      \\
    \hline
    \end{tabular}
  \caption{Vectorizations in different TDA packages}
  \label{tab:vects}

\end{table}

The TDArea R package and its Python counterpart aim to fill these gaps. Its contributions can be summarized as following:
\begin{enumerate}
\item It expands the list of implemented vector summaries of PDs by providing vectorization of eight functional summaries found in the TDA literature: \emph{Betti function}, \emph{persistence landscape function}, \emph{persistence silhouette function}, \emph{persistent entropy summary function} \cite{atienza2020stability}, \emph{Euler characteristic curve} \cite{richardson2014efficient}, \emph{normalized life curve} \cite{chung2022persistence}, \emph{persistence surface} \cite{adams2017persistence} and \emph{persistence block} \cite{chan2022computationally}. In addition, in the updated version o the article such vectorizations as \emph{simple statistics}, \emph{Tropical Coordinates}, and \emph{Template Functions} \cite{TentFunction} were added. In the Appendix you an see definitions of these vectorizations.
\item A univariate summary function $f$ of a PD is typically vectorized by evaluating it at each point of a superimposed one dimensional grid and arranging the resulting values into a vector:
  \begin{equation}\label{stand_vec}
		(f(t_1),f(t_2),\ldots,f(t_n))\in {R}^n,
\end{equation}
where $t_1,t_2,\ldots,t_n$ form an increasing sequence of scale values. For example, the \texttt{landscape()} and \texttt{silhouette()} functions of the \texttt{TDA} package compute vector summaries of persistence landscapes and silhouettes in this manner. The \texttt{TDAvec} package instead employs a different vectorization scheme which involves computing the average values of $f$ between two consecutive scale values $t_i$ and $t_{i+1}$ using integration: 
\begin{equation} 
	\Big(\frac{1}{\Delta t_1}\int_{t_1}^{t_2}f(t)dt,\frac{1}{\Delta t_2}\int_{t_2}^{t_3}f(t)dt,\ldots,\frac{1}{\Delta t_{n-1}}\int_{t_{n-1}}^{t_n}f(t)dt\Big)\in{R}^{n-1}, 
\end{equation}
where $\Delta t_i=t_{i+1}-t_i$. Unlike (\ref{stand_vec}), this vectorization method does not miss the behavior of $f$ between neighboring scale points and applies to all univariate summary functions which are easy to integrate, namely persistence silhouette, persistent entropy summary function, Euler characteristic curve, normalized life curve and Betti function. 
\item To achieve higher computational efficiency, all code behind the vector summaries of \texttt{TDAvec} is written in C++. For example, in \texttt{R}, computing the persistence landscape from a PD with the \texttt{TDAvec} package is more than 200 times faster than with the \texttt{TDA} package.
\item  We created a Python library implementing all TDAvec functions, thus extending their availability to the Python ecosystem. The library offers the exact same interface, and a scikit-learn-type interface is also provided for completeness.
\end{enumerate}

The \texttt{TDAvec} \texttt{R} package and a vignette showing its basic usage with examples are available on the CRAN repository\footnote{https://cran.r-project.org/web/packages/TDAvec/index.html}. For \texttt{Python} examples, we refer the readers to sample notebook presented in \cite{pyTDAvec:2024}.

\subsection{Usage Examples}

Let us first describe how the R library TDAvec can be installed and used.

The current version of this library is available on CRAN, so the simplest way to install it is by using the standard R method:

\begin{verbatim}
   > install.packages("TDAvec") 
\end{verbatim}

Alternatively you can install it from the GitHub repository:

\begin{verbatim}
    > devtools::install_github("uislambekov/TDAvec")
\end{verbatim}

After downloading the library, you can use functions such as computePersistenceLandscape, computePersistenceSilhouette, etc., to calculate the corresponding vectorization summaries of persistence diagrams.

\begin{figure}
  \centering
  \includegraphics[width = 0.45\textwidth]{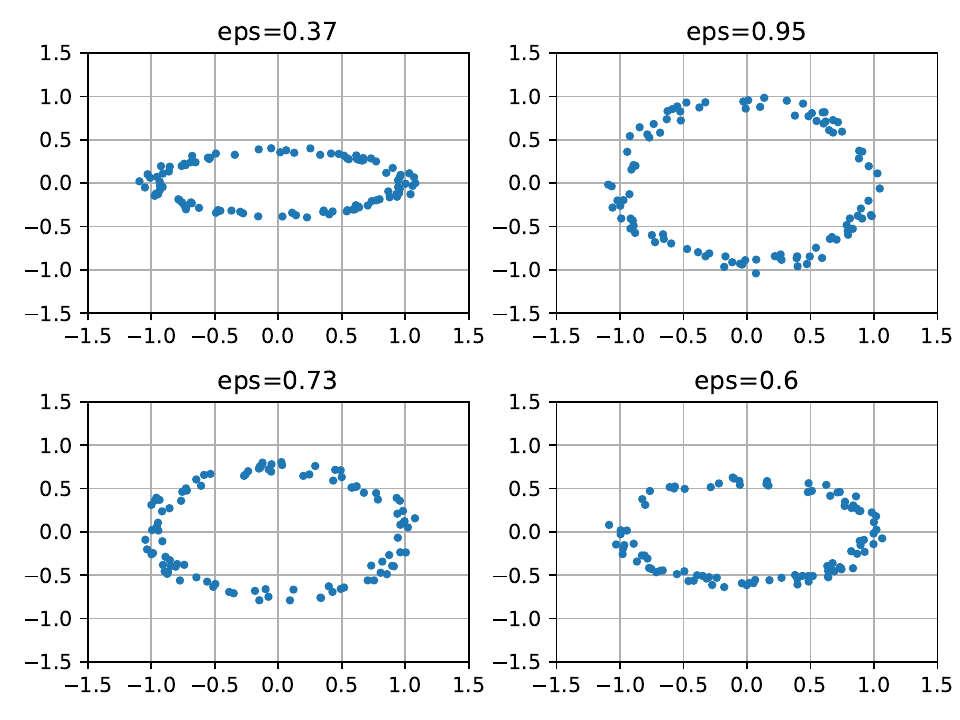}%
  \includegraphics[width = 0.45\textwidth]{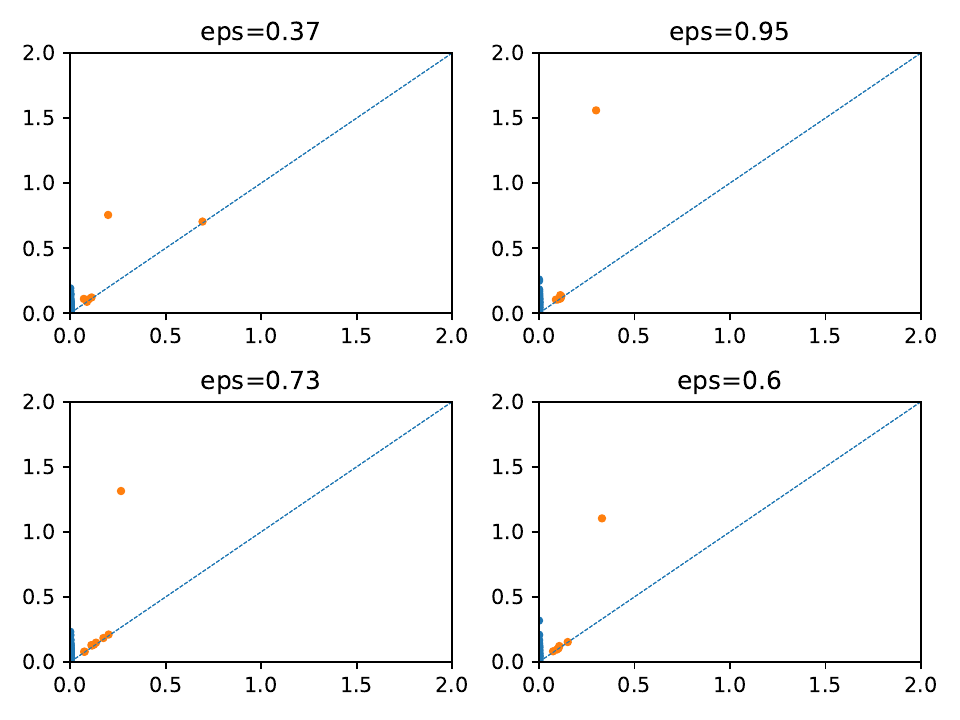}
  \caption{Examples o the point clouds and PDs}
  \label{fig:XPDs}
\end{figure}

Suppose we have a random set of squeeze factors $e^a \in [0,1]$, and for each of them, we create a cloud of points located around an ellipse:
$$
(x_i, y_i)^a = ( r_i \cos\phi_i, e^a r_i\sin\phi_i)
$$
Here are R commands to reate the point clouds:
\begin{verbatim}
    > epsList <- round(as.vector(read.csv("./epsList.csv", header = FALSE))[[1]], 3)
    > clouds <- lapply(epsList, function(e) createEllipse(100, a=1, b=e))
    > PDs <- lapply(1:length(clouds), function(i) 
           ripsDiag(clouds[[i]], maxdimension = 1, maxscale = 2)$diagram
        )
\end{verbatim}
In the figure \ref{fig:XPDs} you can see examples of the created data and corresponding persistence diagrams.

\begin{figure}
  \centering
  \includegraphics[width = 0.5\textwidth]{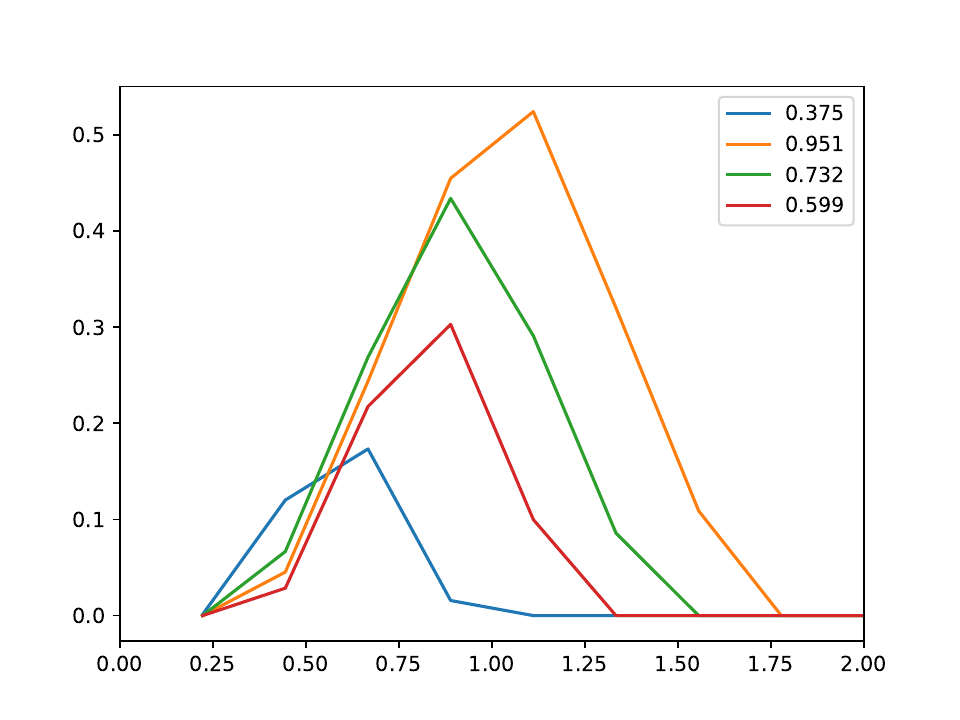}
  \caption{Persistence landscape summaries created from presented above point clouds}
  \label{fig:PLs}
\end{figure}

Created point clouds can be converted into persistence diagrams using such functions as \texttt{ripsDiag} from \texttt{TDA} package. Each of the for each of the diagrams we can calculate Persistence Landscape summary with the help of \texttt{computePerststenceLandscape(diagram, homDim, x)} function\footnote{Note that the names of the vectorization functions have changed during the last update}. In figure \ref{fig:PLs} you can see the result.

The python interface of \texttt{TDAvec} library can be installed from the GitHub repository

\begin{verbatim}
    > pip install git+https://github.com/ALuchinsky/tdavec
\end{verbatim}
Note that it was tested on Python 3.11 and uses numpy (version numpy==1.2) and ripser (version 0.6.8) libraries.

You can directly access all available in R interface unctions of the package, like \texttt{computePerststenceLandscape}, etc. In addition, python implementation follows the scikit-learn interface, so it uses a class-oriented approach. After installing the library, you can load all the objects with:

\begin{verbatim}
    > from tdavec.TDAvectorizer import TDAvectorizer, createEllipse
\end{verbatim}
You can now use an object of the class TDAvectorizer to create persistence diagrams and vectorize them. Here is an example of Python code that performs the same job as the R commands above:

\begin{verbatim}
    > epsList = np.random.uniform(low = 0, high = 1, size = 500)
    > clouds = [createEllipse(a=1, b=eps, n=100) for eps in epsList]
    > v = TDAvectorizer()
    > v.fit(clouds)
\end{verbatim}
The last line effectively uses the ripser library to convert data point clouds into a list of persistence diagrams, and the results will be saved in the v.diags property of the vectorizer object. To perform vectorization, we can simply call the fit() method of the same object:

\begin{verbatim}
    > v.setParams({"scale":np.linspace(0, 2, 10)})
    > X = v.transform(output="PL", homDim=1)
\end{verbatim}
Note that the scale sequence (i.e., grid points for persistence landscape calculation) is specified in the first line. The type of vectorization summary can be specified as a parameter. This approach makes it extremely convenient to easily switch between methods in the analysis and assess accuracy. Here's how we can conduct such an analysis.

\begin{figure}
  \centering
  \includegraphics[width = 0.7\textwidth]{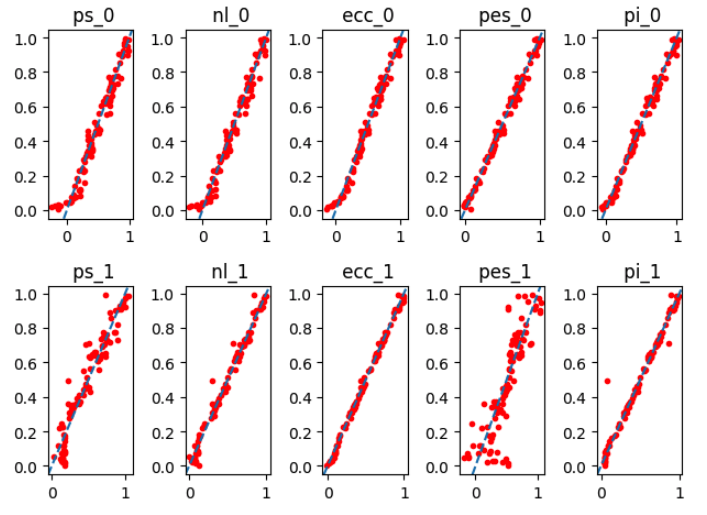}%
  \caption{Linear Regression comparison}
  \label{fig:LR_comparison}
\end{figure}

First, we define a helper function that vectorizes persistence diagrams using the specified method, builds a linear regression model based on the obtained predictors, and reports accuracies on the training and testing subsets:

\begin{verbatim}
    > def makeSim(method, homDim, vec = v, y=epsList):
    >    X =v.transform(output=method, homDim=homDim)
    >    Xtrain, Xtest, ytrain, ytest = train_test_split(X, y, 
                    train_size=0.8, random_state=42)
    >    model = LinearRegression().fit(Xtrain, ytrain)
    >    test_preds = model.predict(Xtest)
    >    score = model.score(Xtest, ytest)
    >    res = {"method":method, "homDim":homDim, 
          "test_preds":test_preds, "y_test":ytest, "score":score}
    >    return res
\end{verbatim}
Now we can use this function to collect data, construct accuracy tables, and compare true values and model predictions for each vectorization method:

\begin{verbatim}
    > methodList = ["vab", "ps", "nl", "ecc", "fda"]
    > df = pd.DataFrame()
    > for homDim in [0, 1]:
    >    for method in methodList:
    >    df = pd.concat([df, pd.DataFrame( makeSim(method, homDim))])
\end{verbatim}
Using the created data frame, we can easily compare performance accuracies (see table \ref{tab:acc}). It's also simple to check the accuracy of the models by inspecting true vs. predicted scatter plots (see figure \ref{fig:LR_comparison}).

\begin{table}[htbp]
  \centering
  \begin{tabular}{lrr}
   \hline
   Method & dim 0 & dim 1 \\
   \hline
   ECC & 0.976 & 0.996 \\
   FDA & 0.983 & 0.985 \\
   NL & 0.8930 & 0.980 \\
   PS & 0.903 & 0.914 \\
   VAB & 0.976 & 0.985 \\
   \hline
  \end{tabular}
  \caption{Test Accuracies of different ML methods}
  \label{tab:acc}
\end{table}

\section{Conclusion}

In this paper, we have presented a new software package aimed at facilitating the vectorization of persistence diagrams, a crucial step in leveraging persistent homology for machine learning applications. By addressing the challenges posed by the non-Hilbert nature of persistence diagrams, our package offers an accessible and efficient framework for transforming these topological summaries into forms compatible with standard machine learning workflows. Through practical examples, we have demonstrated the utility and flexibility of the package, highlighting its potential to accelerate research and application in topological data analysis. We believe this tool will serve as a valuable resource for practitioners and researchers alike, enabling broader adoption of TDA techniques in various domains.

Future work could explore extending the package to integrate novel vectorization techniques and support a wider range of applications, further advancing the impact of TDA in data science and beyond.

The authors would like to thank Dr. Kit Chan and Dr. Rebeca Sanders for fruitful discussions.

\bibliographystyle{alpha}
\bibliography{paper}

\renewcommand{\etalchar}[1]{$^{#1}$}

\appendix
\section{Definitions of the summary functions in \texttt{TDAvec}}

Let $D=\{(b_i,d_i)\}_{i=1}^N$ be a persistence diagram. 

1) The $k$-th order \textbf{persistence landscape} function of $D$ is defined as $$\lambda_k(t) = k\hbox{max}_{1\leq i \leq N} \Lambda_i(t), \quad k\in N,$$
where $k\hbox{max}$ returns the $k$th largest value and 
$$\Lambda_i(t) = \left\{
        \begin{array}{ll}
            t-b_i & \quad t\in [b_i,\frac{b_i+d_i}{2}] \\
            d_i-t & \quad t\in (\frac{b_i+d_i}{2},d_i]\\
            0 & \quad \hbox{otherwise}
        \end{array}
    \right.$$

2) The $p$th power \textbf{persistence silhouette} function:
$$\phi_p(t) = \frac{\sum_{i=1}^N |d_i-b_i|^p\Lambda_i(t)}{\sum_{i=1}^N |d_i-b_i|^p},$$
where
$$\Lambda_i(t) = \left\{
        \begin{array}{ll}
            t-b_i & \quad t\in [b_i,\frac{b_i+d_i}{2}] \\
            d_i-t & \quad t\in (\frac{b_i+d_i}{2},d_i]\\
            0 & \quad \hbox{otherwise}
        \end{array}
    \right.$$

3) The \textbf{persistent entropy summary} function is defined in the following way:
$$
S(t)=-\sum_{i=1}^N \frac{l_i}{L}\log_2{(\frac{l_i}{L}})\mathbf 1_{[b_i,d_i]}(t),
$$ where $l_i=d_i-b_i$ and $L=\sum_{i=1}^Nl_i$.

4) The \textbf{Betti Curve} vectorization: 
$$
\beta(t)=\sum_{i=1}^N w(b_i,d_i)\mathbf 1_{[b_i,d_i)}(t),
$$ where the weight function $w(b,d)\equiv 1$.

5) The \textbf{Euler characteristic curve}: 
$$
\chi(t)=\sum_{k=0}^d (-1)^k\beta_k(t),
$$ where $\beta_0,\beta_1,\ldots,\beta_d$ are the Betti curves corresponding to persistence diagrams $D_0,D_1,\ldots,D_d$ of dimensions $0,1,\ldots,d$ respectively, all computed from the same filtration. Note that unlike all other vectorizations the Euler characteristics curve is calculated using all available dimensions.

6) The \textbf{normalized life curve}: 
$$
sl(t)=\sum_{i=1}^N \frac{d_i-b_i}{L}\mathbf{1}_{[b_i,d_i)}(t),
$$ where $L=\sum_{i=1}^N (d_i-b_i)$.

7) \textbf{Persistence Surface}:
$$\rho(x,y)=\sum_{i=1}^N f(b_i,p_i)\phi_{(b_i,p_i)}(x,y),$$ where $\phi_{(b_i,p_i)}(x,y)$ is 
the Gaussian distribution with mean $(b_i,p_i)$ and 
covariance matrix $\sigma^2 I_{2\times 2}$ and 
$$
f(b,p) = w(p)=\left\{
        \begin{array}{ll}
            0 & \quad p\leq 0 \\
            p/p_{max} & \quad 0<p<p_{max}\\
            1& \quad p\geq p_{max}
        \end{array}
    \right.
$$
is the weighting function with $p_{max}$ being the maximum persistence value among all persistence diagrams considered in the experiment.

8) The \textbf{persistence block}:

$$
f(x,y)=\sum_{i=1}^N \mathbf 1_{E(b_i,p_i)}(x,y),
$$
where $E(b_i,p_i)=[b_i-\frac{\lambda_i}{2},b_i+\frac{\lambda_i}{2}]\times [p_i-\frac{\lambda_i}{2},p_i+\frac{\lambda_i}{2}]$ and $\lambda_i=2\tau p_i$ with $\tau\in (0,1]$.

9) \textbf{Simple statistics}:

\newcommand{\eqn}[1]{$#1$}
\newcommand{\code}[1]{\texttt{#1}}

For a given persistence diagram \code{computeStats()} calculates descriptive statistics of the birth, death, midpoint (the average of birth and death), and lifespan (death minus birth) values corresponding to a specified homological dimension. Additionally, it computes the total number of points and entropy of the lifespan values. Points in \eqn{D} with infinite death values are ignored.

10) \textbf{Tropical Coordinates}:

This vectorization for each homological dimension  computes the following seven tropical coordinates based on the lifespans (or persistence) \eqn{\lambda_i = d_i - b_i}:

\itemize{
  \item \eqn{F_1 = \max_i \lambda_i}.
  \item \eqn{F_2 = \max_{i<j} (\lambda_i+\lambda_j)}.
  \item \eqn{F_3 = \max_{i<j<k} (\lambda_i+\lambda_j+\lambda_k)}.
  \item \eqn{F_4 = \max_{i<j<k<l} (\lambda_i+\lambda_j+\lambda_k+\lambda_l)}.
  \item \eqn{F_5 = \sum_i \lambda_i}.
  \item \eqn{F_6 = \sum_i \min(r \lambda_i, b_i)}, where \eqn{r} is a positive integer.
  \item \eqn{F_7 = \sum_j \big(\max_i(\min(r \lambda_i, b_i)+\lambda_i)  - (\min(r \lambda_j, b_j)+\lambda_j)\big)}.
}

11) \textbf{Template Function Vectorization}

This vectorization uses a collection of tent template functions defined by

\begin{align*}
  G_{(b,p),\delta}(D) &=
                        \sum_{i=1}^N\max \left\{ 0, 1 - \frac{1}{\delta} \max \left( |b_i - b|, |p_i - p| \right) \right\},  
\end{align*}
where \eqn{p_i=d_i-b_i} (persistence), \eqn{b\geq 0}, \eqn{p>0} and  \eqn{0<\delta<p}. The point \eqn{(b,p)} is referred to as the center. Points in \eqn{D} with infinite death values are ignored.

12) \textbf{Algebraic Functions}

This vectorization computes following four algebraic functions bases on birth and death values:
\begin{align}
  f_1 &= \sum_i b_i\left(d_i - b_i\right) \\
  f_2 &= \sum_i \left( d_{\max} - d_i\right) \left(d_i - b_i\right) \\
  f_3 & = \sum_i b_i^2\left(d_i - b_i\right)^4 \\
  f_4 &= \sum_i \left( d_{\max} - d_i\right)^2 \left(d_i - b_i\right)^4,
\end{align}
where $d_{\max} = \max d_i$. Points in $D$ with infinite death value are ignored.

13) Complex Polynomial

This vectorization computes the coefficients of a complex polynomial
\begin{align}
  C_X(z) = \prod_i \left[z - X\left(b_i, d_i\right)\right],
\end{align}
where $X: R \to C$ is any of the following three functions:
\begin{align}
  R(x, y) &= x + i y,\\
  S(x, y) &=
            \begin{cases}
              \frac{y-x}{\sqrt{2(x^2+y^2)}} (x+i y), & \textrm{if } (x,y) \ne (0,0) \\
              0 & \textrm{otherwise}
            \end{cases}
\end{align}
Points in $D$ with infinite death value are ignored.

15) \textbf{Functional Data Analysis}:

The drawback of all listed above vectorization methods is that we have to choose some grid points on t axis or on $(b,p)$ plane to get the numeric array of predictors from the vectorized functions. An alternative approach is to select some set of basis functions $\xi_i(t)$ and expand the result of the vectorization (for example, Betti Curve) as a combination of them:
$$
g(t) = \sum_i I_{[b_i,d_i)]}(t) = \sum_k c_k \xi_k(t)
$$
The coefficients $\{c_k\}$ will be predictors in this approach. It is clear for this type of selection we do not need to select arbitrary the grid points. In addition, it is more robust to small changes of the function $g(t)$ or, alternatively, to birth-death points in the persistence diagram.

\end{document}